\documentclass[11pt]{article}
\usepackage{bbm}
\usepackage{graphicx}
\begin{document}

\title{\bf On Monge-Kantorovich Problem in the Plane}
\author{\textsc{{\small By \ Yinfang SHEN \ and \ Weian ZHENG\thanks{%
Research partially supported N.S.F.C. Grant 10771070}}}}
\date{}
\maketitle
\renewcommand{\theequation}{\arabic{section}.\arabic{equation}}
\newtheorem{theo}{Theorem}[section]
\newtheorem{defi}[theo]{Definition}
\newtheorem{lem}[theo]{Lemma}
\newtheorem{prop}[theo]{Proposition}
\newtheorem{coro}[theo]{Corollary}
\newtheorem{rem}[theo]{Remark}
\newtheorem{cond}[theo]{Condition}
\newtheorem{exam}[theo]{Example}
\newcounter{bean}

\begin{abstract} We transfer the celebrating Monge-Kontorovich problem
in a bounded domain of Euclidean plane into a Dirichlet boundary
problem associated to a quasi-linear elliptic equation with
$0-$order term missing in its diffusion coefficients:
\begin{eqnarray*}
A(x,\ F'_x)F''_{xx}+B(y,\ F'_y)F''_{yy}&=&C(x,\ y,\ F'_x,\ F'_y)
\end{eqnarray*}
where $A(.,.)>0, B(.,.)>0$ and $C$ are functions based on the
initial distributions, $F$ is an unknown probability distribution
function and therefore closed the former problem.
\end{abstract}

\vskip .3in

The mass transport problem was first formulated by Monge in 1781,
and concerned finding the optimal way, in the sense of minimal
transportation cost of moving a pile of soil from one site to
another. This problem was given a modern formulation in the work of
Kantorovitch and so is now known as the Monge-Kontorovich problem.

This type of problem has appeared in economics, automatic control,
transportation, fluid dynamics, statistical physics, shape
optimization, expert system, meteorology and financial mathematics.
For example, for the general tracking problem, a robust and reliable
object and shape recognition system is of major importance. A key
way to carry this out is via template matching. which is the
matching of some some object to another within a given catalogue of
objects. Typically, the match will not be exact and hence some
criterion is necessary to measure the \lq\lq goodness of fit\rq\rq .

Many mathematicians from different fields are interested in
Monge-Kontorovich problem. This classical problem was revived in the
mid eighties by the work of Y.Brenier(\cite{B87}, \cite{B91}), who
characterized the optimal transfer plans in terms of gradients of
convex functions. In the last decades, this problem has been
recovered to have a close relationship with certain evolutionary
PDE's, which can be interpreted as gradient flows of certain entropy
functionals with respect to a metric (which is well-known to
probabilists, see Knott-Smitt \cite{KS87} and Rachev-Rueschendorf
\cite{RR98}) involving optimal transportation called Wasserstein
metric. The first application to mathematical physics (kinetic
models) is due to Tanaka, in the seventies. In the early nineties
the use of entropy functionals as a tool to prove convergence to
equilibrium received a strong impulse due to the work of Cercignani,
Carlen, Carvalho, Pulvirenti, Desvillettes, Toscani, Villani and
others. Moreover, Toscani proved that similar methods could be used
to prove optimal convergence to similarity for diffusion equations.
At the same time Jordan, Kinderlehrer and Otto \cite{JKO98}
discovered that the Fokker Planck equation can be solved by a
steepest descent method involving a logarithmic entropy functional
and the Wasserstein distance. This work marks the beginning of the
modern gradient flow theory on Wasserstein spaces. After a few
years, Arnold, Carrillo, Del Pino, Dolbeault, J\"ungel, Markowich,
Toscani and Unterreiter established the link between convergence to
equilibrium for linear and nonlinear Fokker-Planck type equations
and logarithmic Sobolev inequalities, by developing a previous idea
of Bakry and Emery \cite{BE85} (with applications to the Porous
medium equation). The key ingredients of this theory (the
log-Sobolev and the Csiszar-Kullback inequalities) are related to
certain Gaussian isoperimetric inequalities (see e.g. Talagrand and
Otto-Villani). Otto realized simultaneously that nonlinear diffusion
equations can be seen as gradient flows in the 2-Wasserstein space
of probability measures of a free energy functional. This metric
structure has been made rigorous by Ambrosio-Gigli-Savar¨¨. At the
same time, Carrillo-McCann-Villani applied these ideas to granular
media models producing these arguments in smooth settings. A basic
ingredient of this theory is the notion of convexity along geodesics
in the Wasserstein space introduced by McCann, also called
¡°displacement convexity¡±.

Another striking application of the optimal transportation (from the
probabilistic point of view based on martingales theory) is the
justification of the mean field limits of certain stochastic
particle models by means of the theory of concentration inequalities
developed (among the others) by L\'evy, Gromov, Milman, Bobkov,
Ledoux, Malrieu. A computational method for finding entropy
functionals for evolutionary equations has been recently proposed by
Juengel-Mattes. The use of the Wasserstein distance has been also
extended to scalar conservation laws (Bolley-Brenier-Loeper,
Carrillo-Di Francesco-Lattanzio). The use of these ideas to study
the long-time asymptotics of dissipative homogeneous kinetic models
is based on the almost equivalence of the Euclidean transportation
metric with Fourier-based metrics and on the basic mechanism of
contraction of probability metrics (Gabetta-Toscani-Wennberg,
Bisi-Carrillo-Toscani, Bolley-Carrillo, Carrillo-Toscani). Several
(important) authors have been involved in literature of the optimal
transportation theory, with remarkable applications, we mention here
Caffarelli \cite{C92}\cite{C96}\cite{C02}, Ledoux \cite{L01}, Evans
and Gangbo \cite{EG99}, Carlen and Gangbo \cite{CG03}. We recommend
Caffarelli's address to ICM2002 \cite{C02} and Trudinger's invited
lecture to ICM2006 \cite{T06} and L.C.Evans and W.Gangbo's paper
\cite{EG99} for major references from PDE point of view. We also
recommend S.T.Rachev and L.R\"uschendorf's book \cite{RR98} for a
major reference from probability point of view.

\vskip .4in

There are several formulation of Monge-Kontorovich problem. We are
going to use its formulation in terms of probability theory, which
is to find the Kantorovich-Rubinstein-Wasserstein distance in the
plane. Suppose that we are given two probability distributions $P$
and $\tilde P$ on $R^2$. A $4-$dimensional random vector $(X,\
\tilde X)$ with $P$ and $\tilde P$ as the marginal distributions is
called a coupling of this pair $(P,\ \tilde P).$ The minimum of the
coupling distance $\| X-\tilde X\|_{L_2} $ among all such possible
couplings is called Kantorovich-Rubinstein-Wasserstein
$L_2-$distance between $P$ and $\tilde P.$ From weak convergence
theory, it is easy to see the existence of this optimal coupling
$(X,\ \tilde X).$ The problem is to find a concrete way to get them.
It has important applications in both probability theory and mass
transfer problems. However, the problem has been only completely
solved in one dimensional case. In $R^1,$
Kantorovich-Rubinstein-Wasserstein $L_2-$distance is just given by
\cite{va}
\begin{eqnarray}\label{vaa}
\sqrt{\int_0^1|F^{(-1)}(t)-\tilde F^{(-1)}(t)|^2dt}
\end{eqnarray}
where $F$ and $\tilde F$ are distribution functions of $P$ and
$\tilde P$ respectively, $F^{-1}(t)$ and $\tilde F^{-1}(t),$ $(0\leq
t\leq 1)$ are their right inverses.

\vskip .3in

Without losing generality, we may just consider two probability
measures $P$ and $Q$ on $[0,1]\times [0,1].$ Let $X$ and $Y$ be two
random vectors defined on a same probability space with $P$ and
$\tilde P$ as their individual laws. Denote
\[\tilde X=(\tilde X_1, \tilde X_2)=(Y_1+1,Y_2+1)\]
and denote by $\tilde P$ its probability distribution which is on
$[1,2]\times [1,2].$ Then
\begin{eqnarray*}
& &E[|X_1-\tilde X_1|^2+|X_2-\tilde
X_2|^2]\nonumber\\
&=&E[|X_1-Y_1-1|^2+|X_2-Y_2-1|^2]\nonumber\\
&=&E[|X_1-Y_1|^2+|X_2-Y_2|^2]-2E[X_1]+2E[Y_1]-2E[X_2]+2E[Y_2]+2
\end{eqnarray*}
which gives the relation between Kantorovich-Rubinstein-Wasserstein
$L_2-$distance of $(P,\ Q)$ and that of $(P,\ \tilde P).$ Since
$-E[X_1]+E[Y_1]-E[X_2]+E[Y_2]$ is given, it is sufficient to discuss
the later.

Assume that $P$ is a probability measures on $[0,1]\times [0,1]$ and
and $\tilde P$ is a probability measure on $[1,2]\times [1,2]$.
Suppose that the couple $X=(X_1, X_2)$ and $\tilde X=(\tilde X_1,
\tilde X_2)$ give the desired Kantorovich-Rubinstein-Wasserstein
$L_2-$distance. If we denote $Z=(X_1, \tilde X_2),$ then
\[ E[|X-\tilde X|^2]=E[|X-Z|^2]+E[|Z-\tilde X|^2].\]

So it is sufficient to find the distribution of $Z$ which is
supported in $[0,1]\times[1,2].$ We assume further the density
functions $f(x,y)$ of $X$ and $\tilde f(x,y)$ of $\tilde X$ are
smooth and strictly positive on their domains. Denote the marginal
densities
\[ f_1(x)=\int_0^1f(x,y)dy,\hskip .5in f_2(y)=\int_0^1f(x,y)dx\]
and
\[ \tilde f_1(x)=\int_1^2\tilde f(x,y)dy,\hskip .5in \tilde f_2(y)=\int_1^2\tilde f(x,y)dx.\]

Furthermore, denote the conditional distributions
\[ F_1(x|y)={1\over f_2(y)}\int_0^x f(u,y)du,\hskip .5in F_2(y|x)
={1\over f_1(x)}\int_1^y f(x,u)du,\] and
\[ \tilde F_1(x|y)={1\over \tilde f_2(y)}\int_0^x \tilde f(u,y)du,\hskip .5in \tilde F_2(y|x)={1\over \tilde f_1(x)}\int_1^y
\tilde f(x,u)du,\] which are strictly increasing with respect to
their first argument so their inverse functions with respect to
their first arguments exist and denoted as
$G(1,s,y)=F_1^{(-1)}(s|y),$ $G(2,x,t)=F_2^{(-1)}(t|x),$ $\tilde
G(1,s,y)=\tilde F_1^{(-1)}(s|y)$ and $\tilde G(2,x,t)=\tilde
F_2^{(-1)}(t|x).$ Without losing generality, we assume that there is
a positive constant $c>0$ such that
\begin{eqnarray}\label{uni}
{1\over \tilde f_2(.)}\tilde G_x(1,.,.)>c\hskip .2in {1\over
f_1(.)}G_y(2,.,.)>c
\end{eqnarray}
and that all functions appeared in the later equation (\ref{hh2})
are sufficiently smooth. Our those regularity hypotheses will not
affect the generality of our problem, because what we will treat
later is the unknown distribution function $F$ of $Z,$ which is
continuous under the weak convergence of the laws of $(X,\tilde X).$
Therefore we can always use the usual regularizing approximation
procedures.

Denote by ${\cal H}$ the set of all density functions $q(x,y)$ on
$[0, 1]\times [1, 2]$ such that $f_1(x)=\int_1^2q(x,y)dy$ and
$\tilde f_2(y)=\int_0^1q(x,y)dx.$

We are looking for a density function $p(x,y)$  satisfying

1) $p\in {\cal H};$

2) $p(x,y)=q(x,y)$ minimizes
\begin{eqnarray}
& &\int_1^2\int_0^1|s-\tilde G(1,\int_0^s{q(u,y)\over \tilde
f_2(y)}du,
y)|^2q(s,y) ds\, dy\nonumber\\
& &+\int_0^1\int_1^2|t-G(2,x,\int_0^t{q(x,v)\over
f_1(x)}dv)|^2q(x,t) dt\, dx
\end{eqnarray}

\vskip .4in

For $0<a<a_1<1$ and $1<b<b_1<2$ when $\epsilon $ is small enough,
\[ a+\epsilon <a_1<a_1+\epsilon <1,\hskip .3in b+\epsilon <b_1<b_1+\epsilon
<2.\] Define
\begin{eqnarray} \xi (s,t)=& &I_{([a, a+\epsilon
]\times [b, b+\epsilon ])\cup ([a_1, a_1+\epsilon ]\times [b_1,
b_1+\epsilon
])}(s,t)\nonumber\\
& &-I_{([a, a+\epsilon ]\times [b_1, b_1+\epsilon ])\cup ([a_1,
a_1+\epsilon ]\times [b, b+\epsilon ])}(s,t).
\end{eqnarray}
Then $p(s,t)+\delta \xi (s,t)\in {\cal H}$ when both $\epsilon ,
\delta $ are small. Since $p$ is the minimum, by (\ref{vaa})
\begin{eqnarray}\label{uu}
0\leq & &{1\over \epsilon^2}\int_1^2\int_0^1|s-\tilde
G(1,\int_0^s{p(u,y)+\delta \xi (u,y)\over \tilde f_2(y)}du,
y)|^2(p(s,y)+\delta \xi (s,y)) ds\, dy\nonumber\\
& &+{1\over
\epsilon^2}\int_0^1\int_1^2|t-G(2,x,\int_0^t{p(x,v)+\delta \xi
(x,v)\over
f_1(x)}dv)|^2(p(x,t)+\delta \xi (x,t)) dt\, dx\nonumber\\
& &-{1\over \epsilon^2}\int_1^2\int_0^1|s-\tilde
G(1,\int_0^s{p(u,y)\over \tilde f_2(y)}du,
y)|^2p(s,y) ds\, dy\nonumber\\
& &-{1\over \epsilon^2}\int_0^1\int_1^2|t-G(2,x,\int_0^t{p(x,v)\over
f_1(x)}dv)|^2p(x,t) dt\, dx\\
=& &{1\over \epsilon^2}\{ \int_1^2\int_0^1|s-\tilde
G(1,\int_0^s{p(u,y)+\delta \xi (u,y)\over \tilde f_2(y)}du,
y)|^2p(s,y) ds\, dy\nonumber\\
& &-\int_1^2\int_0^1|s-\tilde G(1,\int_0^s{p(u,y)\over\tilde
f_2(y)}du,
y)|^2p(s,y) ds\, dy\}\nonumber\\
& &+{1\over
\epsilon^2}\{\int_0^1\int_1^2|t-G(2,x,\int_0^t{p(x,v)+\delta \xi
(x,v)\over
f_1(x)}dv)|^2p(x,t) dt\, dx\nonumber\\
& &-\int_0^1\int_1^2|t-G(2,x,\int_0^t{p(x,v)\over
f_1(x)}dv)|^2p(x,t) dt\, dx\}\nonumber\\
& &+{1\over \epsilon^2}\{\int_1^2\int_0^1|s-\tilde
G(1,\int_0^s{p(u,y)+\delta \xi (u,y)\over \tilde f_2(y)}du,
y)|^2\delta \xi (s,y) ds\, dy\nonumber\\
& &+\int_0^1\int_1^2|t-G(2,x,\int_0^t{p(x,v)+\delta \xi (x,v)\over
f_1(x)}dv)|^2\delta \xi (x,t) dt\, dx\}\nonumber
\end{eqnarray}
Letting $\epsilon\to 0,$ we get
\begin{eqnarray}
0\leq & &-2\int_a^{a_1}(\tilde G(1,\int_0^s{p(u,b_1)\over \tilde
f_2(b_1)}du, y)-s)\tilde G'_x(1,\int_0^s{p(u,b_1)\over \tilde
f_2(b_1)}du, b_1){p(s,b_1)\over \tilde f_2(b_1)} ds
\nonumber\\
& &+2\int_a^{a_1}(\tilde G(1,\int_0^s{p(u,b)\over \tilde f_2(b)}du,
b)-s)\tilde G'_x(1,\int_0^s{p(u,b)\over \tilde f_2(b)}du,
b){p(s,b)\over \tilde f_2(b)} ds \nonumber\\
& &-2\int_b^{b_1}(G(2,a_1,\int_0^t{p(a_1,v)\over
f_1(a_1)}dv)-t)G'_y(2,a_1,\int_0^t{p(a_1,v)\over
f_1(a_1)}dv){p(a_1,t)\over
f_1(a_1)} dt\nonumber\\
& &+2\int_b^{b_1}(G(2,a,\int_0^t{p(a,v)\over
f_1(a)}dv)-t)G'_y(2,a,\int_0^t{p(a,v)\over f_1(a)}dv){p(a,t)\over
f_1(a)} dt\nonumber\\
& &+|a_1-\tilde G(1,\int_0^{a_1}{p(u,b_1)\over \tilde f_2(b_1)}du,
b_1)|^2 -|a-\tilde G(1,\int_0^a{p(u,b_1)\over\tilde  f_2(b_1)}du,
b_1)|^2
\nonumber\\
& &-|a_1-\tilde G(1,\int_0^{a_1}{p(u,b)\over \tilde f_2(b)}du, b)|^2
+|a-\tilde G(1,\int_0^a{p(u,b)\over \tilde f_2(b)}du, b)|^2
\nonumber\\
& &+|b_1-G(2,a_1,\int_0^{b_1}{p(a_1,v)\over f_1(a_1)}dv)|^2
-|b-G(2,a_1,\int_0^b{p(a_1,v) \over f_1(a_1)}dv)|^2\nonumber\\
& &-|b_1-G(2,a,\int_0^{b_1}{p(a,v)\over f_1(a)}dv)|^2
+|b-G(2,a,\int_0^b{p(a,v) \over f_1(a)}dv)|^2\nonumber
\end{eqnarray}
Multiplying both sides by ${1\over (a_1-a)(b_1-b)},$ letting $(a_1-
a)(b_1- b)\rightarrow 0,$ we get
\begin{eqnarray}\label{eu}
{\partial^2\over \partial x\partial y}M(x,y)\geq 0
\end{eqnarray}
where
\begin{eqnarray}\label{5}
& &M(x,y)\nonumber\\
&= &-2\int_0^x(\tilde G(1,\int_0^s{p(u,y)\over \tilde f_2(y)}du,
y)-s)\tilde G'_x(1,\int_0^s{p(u,y)\over \tilde f_2(y)}du,
y){p(s,y)\over \tilde f_2(y)} ds
\nonumber\\
& &-2\int_1^y(G(2,x,\int_1^t{p(x,v)\over
f_1(x)}dv)-t)G'_y(2,x,\int_1^t{p(x,v)\over f_1(x)}dv){p(x,t)\over
f_1(x)} dt\nonumber\\
& &+|x-\tilde G(1,\int_0^x{p(u,y)\over \tilde f_2(y)}du, y)|^2
+|y-G(2,x,\int_1^y{p(x,v)\over f_1(x)}dv)|^2\nonumber\\
&=&\int_0^x2(s-\tilde G(1,\int_0^s{p(u,y)\over \tilde f_2(y)}du, y))ds\nonumber\\
& &+\int_1^y2(t-G(2,x,\int_1^t{p(x,v)\over f_1(x)}dv))dt\nonumber\\
&=&x^2+y^2-2\int_0^x \tilde G(1,\int_0^s{p(u,1)\over \tilde
f_2(1)}du,1)ds -2\int_1^y
G(2,0,\int_1^t{p(0,v)\over f_1(0)}dv)dt\nonumber\\
& &-2\int_1^y\int_0^x \{ [\tilde G(1,\int_0^s{p(u,t)\over \tilde
f_2(t)}du, t)]'_t+[G(2,s,\int_1^t{p(s,v)\over f_1(s)}dv)]'_s\} ds\,
dt
\end{eqnarray}

\vskip .2in On the other hand, if one replace $p+\delta \xi $ by
$p-\delta \xi ,$ the same computation leads
\begin{eqnarray}\label{eud}
{\partial^2\over \partial x\partial y}M(x,y)\leq 0
\end{eqnarray}
Thus we deduce from (\ref{eu}) and (\ref{eud}) that
\begin{eqnarray*}\label{kkk}
{\partial^2\over \partial x\partial y}M(x,y)= 0\hskip .3in (\forall
\ 0<x<1<y<2),
\end{eqnarray*}
or $(\forall \ 0<x<1<y<2)$
\begin{eqnarray}\label{hh}
[\tilde G(1,\int_0^x{p(u,y)\over \tilde f_2(y)}du,
y)]'_y+[G(2,x,\int_1^y{p(x,v)\over f_1(x)}dv)]'_x= 0\hskip .3in
\end{eqnarray}

\vskip .2in

Denote the probability distribution function
$F(x,y)=\int_0^x\int_1^y p(s,t)\, dt\, ds.$ Then
\[ F'_x(x,y)=f_1(x)\int_1^y {p(x,v)\over f_1(x)}dv,\hskip .3in
F'_y(x,y)=f_2(y)\int_0^x {p(u,y)\over f_2(y)} du,\] (\ref{hh})
becomes
\begin{eqnarray}\label{hh1}
[\tilde G(1,{1\over \tilde f_2(y)}F'_y(x,y), y)]'_y+[G(2,x,{1\over
f_1(x)}F'_x(x,y))]'_x= 0
\end{eqnarray}
or
\begin{eqnarray}\label{hh2}
& &\tilde G_x(1,{1\over \tilde f_2(y)}F'_y(x,y), y){1\over \tilde
f_2(y)}F''_{yy}(x,y)+G_y(2,x,{1\over
f_1(x)}F'_x(x,y)){1\over f_1(x)}F''_{xx}(x,y)\nonumber\\
&=&-\tilde G_y(1,{1\over\tilde  f_2(y)}F'_y(x,y), y)-G_x(2,x,{1\over
f_1(x)}F'_x(x,y))\nonumber\\
& &+\tilde G_x(1,{1\over \tilde f_2(y)}F'_y(x,y), y){f_2'(y)\over
f^2_2(y)}F'_y(x,y)\nonumber\\
& &+G_y(2,x,{1\over f_1(x)}F'_x(x,y)){f_1'(x)\over
f_1^2(x)}F'_x(x,y)
\end{eqnarray}
which is a quasi-linear elliptic equation with unknown $F(x,y),$
which satisfies the uniform ellipticity condition (\ref{uni}). Since
its diffusion coefficients only contain the first order partial
derivatives, under minus regularity condition, the solution of the
Dirichlet boundary problem $(\forall x\in [0,1], \forall y\in
[1,2]):$
\[ F(0,y)=0,\ F(x,1)=0, \ F(x,2)=\int_0^xf_1(s)ds, \
F(1,y)=\int_0^yf_2(t)dt\]
has a unique solution (\cite{gt} p.264).
\vskip .3in

Furthermore, if we plug (\ref{hh}) into (\ref{5}), then
\begin{eqnarray}\label{15}
M(x,y) &=&x^2+y^2-2\int_0^x \tilde G(1,\int_0^s{p(u,1)\over
\tilde f_2(1)}du,1)ds\nonumber\\
& & -2\int_1^y G(2,0,\int_1^t{p(0,v)\over f_1(0)}dv)dt
\end{eqnarray}
That is, $M(x,y)$ can be written as a closed-form solution which
depends only the initial values $p(0,.)$ and $p(.,0)$

\vskip 0.3in

Yinfang Shen

\small{Department of Statistics, East China Normal University,
Shanghai, China, 200062}

\small{Department of Mathematics, University of California, Irvine,
CA 92697, USA}

\vskip .2in

Weian Zheng

\small{Department of Statistics, East China Normal University,
Shanghai, China, 200062}

\small{Department of Mathematics, University of California, Irvine,
CA 92697, USA}

email address: wzheng@uci.edu

\enddocument